\input amssym.tex
%%%%%%%%%%%%%%%%%%%%%%%
\font\twbf=cmbx12

\font\sc=cmcsc10

\def\rond{{\scriptstyle\circ}}

\def\cf{{\it cf.\/}\ }    
\def\ie{{\it i.e.\/}\ }
\def\etc{{\it etc.}} 
\def\up#1{\raise 1ex\hbox{\sevenrm#1}}

\def \bg {\bigskip \goodbreak}

\def\ref#1&#2&#3&#4&#5\par{\par{\leftskip = 5em {\noindent
\kern-5em\vbox{\hrule height0pt depth0pt width
5em\hbox{\bf[\kern2pt#1\unskip\kern2pt]\enspace}}\kern0pt}
{\sc\ignorespaces#2\unskip},\
{\rm\ignorespaces#3\unskip}\
{\sl\ignorespaces#4\unskip\/}\
{\rm\ignorespaces#5\unskip}\par}}

\def \comp{ \Bbb{C} }

\def \adh  {\overline}

\def \abs#1{\left\vert#1\right\vert }

\def\N#1{\muskip0=-2mu{\left|\mkern\muskip0\left|
#1\right|\mkern\muskip0\right|}}

\def\tore{\Bbb{T}}
\def\hot{\widehat\otimes}
\def\wt#1{\widetilde #1}
\def\wcheck#1{\smash{
        \mathop{#1}\limits^{\scriptscriptstyle{\vee}}}}

\def\Otimes_#1^#2{\matrix{{}_{#2}\cr
      \bigotimes \cr {}^{#1}}}

\def\wco#1#2{\lower0.5ex\hbox{$\scriptstyle{{#1}\wcheck\otimes
{#2}}$}}

\def\refl#1&#2&#3&#4&#5\par{\par{\leftskip = 7em {\noindent
\kern-7em\vbox{\hrule height0pt depth0pt width
5em\hbox{\bf[\kern2pt#1\unskip\kern2pt]\enspace}}\kern0pt}
{\sc\ignorespaces#2\unskip},\
{\rm\ignorespaces#3\unskip}\
{\sl\ignorespaces#4\unskip\/}\
{\rm\ignorespaces#5\unskip}\par}}

\def\ref#1&#2&#3&#4&#5\par{
\item{[{\bf\ignorespaces#1\unskip}]}
{\sc\ignorespaces#2\unskip},\
{\rm\ignorespaces#3\unskip}\
{\sl\ignorespaces#4\unskip\/}\
{\rm\ignorespaces#5\unskip}\par}
 
\def\LLongrightarrow{-\mkern-4.8mu-
       \mkern-9.5mu\longrightarrow}
\def\LLongleftarrow{\longleftarrow\mkern-9mu-\mkern-7.0mu-}

\def\mapright#1{\smash{
        \mathop{\LLongrightarrow}\limits^{#1}}}
\def\mapleft#1{\smash{
        \mathop{\LLongleftarrow}\limits^{#1}}}
\def\mapupr#1{\biggl\uparrow
        \rlap{$\vcenter{\hbox{$#1$}}$}}
\def\mapupl#1{\rlap{$\vcenter{\hbox{$#1$}}$}
       \, \biggl\uparrow}

\def\th{\theta}

\def\liminf_#1{\mathop{\underline{\hbox{\rm lim}}}\limits_{#1}}
\def\limsup_#1{\mathop{\overline{\hbox{\rm lim}}}\limits_{#1}}

%%%%%%%%%%%%%%%%%%%%%%%
\pageno=1
{\twbf
\centerline{L'Application Canonique ${\bf  J:\widetilde
H^2(X) \hot \widetilde H^2(X) \rightarrow \widetilde H
^1(X\hot X) }$ }
\centerline{n'est pas  Surjective en G\'en\'eral }
}
\bg
\bg
\bg
\centerline{{\sc Omran Kouba}}\par
\centerline{Department of Mathematics}\par
\centerline{{\sl Higher Institute for Applied Sciences and Technology}}\par
\centerline{P.O. Box 31983, Damascus, Syria.}\par
\centerline{{\it E-mail} : omran\_kouba@hiast.edu.sy}\par
\bg
\bg
{\bf R\'esum\'e :}
On introduit la propri\'et\'e $  H^1$-projective, et on
l'utilise pour construire un espace de Banach $ X $ pour
lequel l'application canonique  
$ J:\wt H^2(X)
\hot \wt H^2(X) \rightarrow \wt H
^1(X\hot X)$  n'est pas surjective.
\par

\bg
{\twbf
\centerline{The Natural\  Map ${\bf  J:\widetilde
H^2(X) \hot \widetilde H^2(X) \rightarrow \widetilde H
^1(X\hot X) }$ }
 \centerline{is not Surjective in General}
}
\bg
\bg
{\bf Abstract :} 
We introduce the $H^1$-projective property, and use it to
construct a Banach space $ X $ such that the natural map
$ J:\wt H^2(X)
\hot \wt H^2(X) \rightarrow \wt H
^1(X\hot X)$ is not onto.

\vskip1cm
\baselineskip=19pt
\noindent{\bf Abridged English Version }
\par
In this Note all Banach spaces considered are complex. Let
$ \tore $ denote the interval $[0,2\pi]$.   For $ p \in
[1,+\infty[ $ and for a banach space $ X $, let $ \wt
H^p(X) $ be  the closed subspace of $ L^p(\tore;X) $
spanned by  analytic polynomials with coefficients in $ X
$.
\par
If $ T: X \rightarrow Y $ is a linear operator between
two Banach spaces. Then the formula
$$
\wt T\left(\sum_{n=0}^m   z^n x_n \right ) = \sum_{n=0}^m
   z^n T(x_n) $$
defines an operator $\wt T: \wt H^p(X) \rightarrow \wt
H^p(Y) $ with the same norm.
\par
We will say that $ X $ is $ H^1$-projective, if there
exists a metric surjection 
$ \sigma : \ell^1(I) \rightarrow X$ --\ie $ \sigma^\ast $
is an isometric embedding--  such that $ \wt \sigma $ is
a surjection from $\wt H^1(\ell^1(I)) $ onto $ \wt H^1(X)
$.  It is easy to see  that $ X $ is $ H^1$-projective,  if
and only if,  $ H^1(\comp)\hot X = \wt H^1(X) $;  and that such
a space is of cotype 2. See [{\bf K }].
\par
For $ F = \sum_0^m h_n \otimes g_n \in \wt H^2(X)
\hot \wt H^2(Y) $, we denote by $ J(F) $ the element in
$\wt H^1(X\hot Y) $ defined by 
$$
J(F)(z) = \sum_0^m h_n(z) \otimes g_n(z) . $$
It is easy to check that
$$
\N{J(F)}_{\wt H^1(X\hot Y)} \leq \N{F}_{ \wt H^2(X)
\hot \wt H^2(Y)} $$
hence one can extend $ J $ to an operator -- still called
$ J $ --  from $\wt H^2(X) \hot \wt H^2(Y) $ into 
$ \wt H^1(X\hot Y) $ of norm one.
\par
Gilles Pisier shows in [{\bf P1}] that this operator is
onto for many couples of Banach spaces $ (X,Y) $,  for
instance if $ X $ and $Y $ are type 2 spaces, or 2-convex
Banach lattices.\par
Our purpose in this Note, is to construct a Banach space $X $ such that the operator $J$ associated to the couple
$ (X,X) $ is not onto. To this end, we use the fact that
if $ X $ is  $H^1$-projective, then $J:\wt H^2(X)
\hot \wt H^2(X) \rightarrow \wt H^1(X\hot X)$ is  onto, if and only if,  $ X\hot X $ is 
$H^1$-projective. It is then enough to find $ X $ which is
$H^1$-projective while $ X\hot X $ is not
$H^1$-projective for some strong reason as being of no
cotype.
The construction is an adaptation of some ideas from 
[{\bf P2}] and [{\bf P3}]. Indeed, we prove the following 
theorem.
\par
\proclaim Theorem. 
Every $ H^1$-projective Banach space $ E $ can be
isometrically embedded in an $ H^1$-projective space $ X
$ such that $ X\hot X = X\wcheck\otimes X $.\par
\bg
It is not difficult to find a space $ X $ satisfying the
preceding theorem and such that both $ X $ and its dual 
$ X^\ast $ are of cotype 2, and  satisfy  Grothendieck's
theorem. Only a minor modification of the construction is
needed to prove this assertion.
\bg

\noindent{\bf 1.D\'{\sevenbf EFINITIONS ET NOTATIONS}:}
\par
Les espaces consid\'er\'es sont des espaces de Banach
complexes.
 Soit $ X $ un espace de Banach, et $ p \in
[1,+\infty[ $, on note $ \wt H^p(X) $ l'adh\'erence dans 
$ L^p(\tore;X) $ des polyn\^omes analytiques \`a
coefficients dans $ X $. ($ \tore $ \'etant $ [0,2\pi]$).
\par
Si $ T:X\rightarrow Y $ est un op\'erateur born\'e d'un
espace de Banach $ X $ dans un autre $ Y $, alors la
formule 
$$\wt T\left(\sum_{n=0}^m  z^n x_n \right ) = \sum_{n=0}^m
 z^n T(x_n) $$
d\'efinit un op\'erateur born\'e $\wt T: \wt H^p(X)
\rightarrow \wt H^p(Y) $ de m\^eme norme.
\par
On dira que $ X $ est $ H^1$-projectif,  s'il existe
un ensemble  $ I $,  une
surjection m\'etrique
$\sigma : \ell^1(I) \rightarrow X $  
(\ie $ \sigma^\ast $ 
est une injection isom\'etrique) et une constante $ K $ 
telle que 
$$
\forall  g \in \wt H^1(X), \quad\exists h \in \wt
H^1(\ell^1(I))\quad :\quad \wt \sigma(h) = g \qquad\hbox{ et }\qquad
\N{h}_{\wt H^1(\ell^1(I)) } \leq K\N{g}_{\wt H^1(X)} 
\eqno (1) $$
pour all\'eger l'\'ecriture, on se contente d'exprimer 
ce qui pr\'ec\`ede en \'ecrivant simplement ``{$ X \hbox{ est
} H_1(K,\sigma,I)$}''.  En effet cette propri\'et\'e est
\'equivalente \`a 
 $ \wt H^1(\comp) \hot X = \wt H^1(X)$,  ce qui justifie
la terminologie ``$X$ est $H^1$-projectif ''.
\par
d'apr\`es une remarque dans [{\bf HP}],  si $ X $ est 
$H^1$-projectif, alors $ X $ v\'erifie aussi (1) en
rempla\c cant $ \wt H^1 $ par $\wt H^p $.
\par
On dira que $(X,Y)$ a la propri\'et\'e $ {\cal P}(c) $ 
si $Y $ est un sous-espace ferm\'e de $ X $ et si 
le  quotient $ Q:X\rightarrow X/Y $ v\'erifie
$$
\forall g \in \wt H^1(X/Y),\quad \exists h \in \wt H^1(X)\quad  :\quad \wt Q(h) = g 
\qquad\hbox{ et }\qquad  \N{h}_{\wt H^1(X)}
\leq c \N{g}_{\wt H^1(X/Y)}  \eqno (2) $$
\par
 Si $ I \subset J $, on notera $ s $ l'injection canonique
 de $\ell^1(I) $ dans  $\ell^1(J) $, d\'efinie par 
$ s(x)(j) = x_j $ si $ j \in I $, et $ s(x)(j) = 0 $ si
 $ j \not\in I $.
\par
Si  $ F = \sum_0^m h_n \otimes g_n \in \wt H^2(X)
\hot \wt H^2(Y) $, on appelle $ J(F) $ l'\'el\'ement de
$\wt H^1(X\hot Y) $ d\'efini par 
$$
J(F)(z) = \sum_{n=0}^m h_n(z) \otimes g_n(z) . $$
Il est facile de v\'erifier que
$$
\N{J(F)}_{\wt H^1(X\hot Y)} \leq \N{F}_{ \wt H^2(X)
\hot \wt H^2(Y)} $$
et donc $ J $ s'\'etend, par densit\'e, en un op\'erateur
(not\'e encore $ J $) de $ \wt H^2(X)
\hot \wt H^2(Y) $ dans $\wt H^1(X\hot Y) $.
\par
Dans [{\bf P1}], Gilles Pisier d\'emontre que cet
op\'erateur est surjectif, pour une large classe de
couples d'espaces $ (X,Y) $, par exemple, si $ X $ et $ Y
$ sont de type 2, des treillis de Banach 2-convexes ou 
des espaces ${\cal L}^\infty$. \par
Le but de cette Note est de donner un exemple d'espace 
$ X $, tel que l'op\'erateur $ J $ associ\'e au couple 
$ (X,X) $ ne soit pas surjectif.
\bg
\noindent{\bf 2.T{\sevenbf H\'EOR\`EMES} :}
\bg
La proposition suivante explique la raison pour laquelle
on a introduit la notion de $ H^1$-projectivit\'e.
\proclaim Proposition 1. 
Soient $ X $ et $ Y $ deux espaces $H^1$-projectifs. Les
assertions suivantes sont \'equivalentes.\hfill\break
\indent 1. $ X \hot Y $ est $ H^1$-projectif.\hfill\break
\indent 2. L'application $ J:\wt H^2(X)
\hot \wt H^2(Y) \rightarrow \wt H
^1(X\hot Y)$  est surjective.
\par
Preuve: 
$1. \Rightarrow 2.$ Ceci est immediat en notant que
$  H^1\hot X \hot Y $ est toujours contenu dans l'image
par $ J $ de $ \wt H^2(X) \hot \wt H^2(Y) $.
\par
$2. \Rightarrow 1.$ Supposons que  $ X \hbox{ est } H_1
(K,\sigma_1,I) $ et que $ Y \hbox{ est } H_1(K,\sigma_2,J)
$. Rappelons que $$\ell^1(I\times J) = \ell^1(I;\ell^1(J)) = 
\ell^1(I)\hot \ell^1(J).$$ Soit $ \sigma : \ell^1(I)\hot 
\ell^1(J)  \rightarrow X \hot Y $ d\'efinie par 
$\sigma(\alpha\otimes \beta) = \sigma_1(\alpha)\otimes
\sigma_2(\beta) $.  $ \sigma $ est clairement une
surjection m\'etrique.
\par
Soit $ F \in \wt H^1(X\hot Y) $, alors d'apr\`es 2. il
existe $ H = \sum_0^\infty h_n \otimes g_n $ dans 
$\wt H^2(X)\hot \wt H^2(Y) $ tel que 
$$ J(H) = F  \qquad\hbox{ et }\qquad     \sum_{n=0}^\infty 
\N{h_n}_{\wt H^2(X)} \N{g_n}_{\wt H^2(Y)}  \leq c \N{F}_
{\wt H^1(X\hot Y) } .$$
D'apr\`es l'hypoth\`ese de $ H^1$-projectivit\'e de $ X $
(resp. $ Y $), et en utilisant la remarque de [{\bf HP}],
 on trouve, pour chaque $ n\geq 0 $ une
fonction $ h_n^\prime \in \wt H^2(\ell^1(I)) $ (resp. 
$ g_n^\prime \in \wt H^2(\ell^1(J)) $), telle que 
$\wt \sigma_1( h_n^\prime)= h_n $ (resp.
$\wt \sigma_2( g_n^\prime)= g_n $), et telle qu'on ait la 
majoration suivante $ \N{h_n^\prime}_{\wt H^2(\ell^1(I))} \leq
 K^\prime \N{h_n}_{\wt H^2(X)}$ (resp. 
 $ \N{g_n^\prime}_{\wt H^2(\ell^1(J))} \leq
 K^\prime \N{g_n}_{\wt H^2(Y)}$).
\par
Consid\'erons l'application canonique  $$ J_1 : \wt
H^2(\ell^1(I)) \hot \wt H^2(\ell^1(J)) \longrightarrow  
\wt H^1(\ell^1(I)\hot\ell^1(J)) ,$$ 
et soit  
$ G = \sum_0^\infty h_n^\prime \otimes g_n^\prime $, on  a
clairement 
$$ \N{ G}_{\wt H^2(\ell^1(I))
\hot \wt H^2(\ell^1(J))} \leq c {K^\prime}^2 \N{F}_
{\wt H^1(X\hot Y) } $$
et
$$\eqalign { \forall  z \in \tore,\qquad\qquad \wt \sigma ( J_1( G))(z)
&= \sigma (J_1(G)(z)) \cr
&=\sum_{n=0}^\infty
\sigma_1(h_n^\prime(z)) \otimes \sigma_2( g_n^\prime(z)) \cr
 &=\sum_{n=0}^\infty h_n(z) \otimes g_n(z) = F(z) } $$
donc $$ \wt \sigma (J_1(G)) = F   \qquad\hbox{ et }\qquad  \N{
J_1(G)}_{\wt H^1(\ell^1(I\times J))} \leq c {K^\prime}^2
\N{F}_ {\wt H^1(X\hot Y) } $$
et ceci d\'emontre la proposition.
\par
\bg
Nous utiliserons aussi le r\'esultat suivant de [{\bf K}]: 
\par
\proclaim Proposition 2. 
Si $ E $ est $ H^1$-projectif, alors $ E $ est de cotype 2.
\par
Il suffit donc de construire  un espace $ H^1$-projectif $ X
$, tel que  $ X\hot X $ ne le soit pas, c'est le cas si ce
dernier n'a pas de cotype. L'espace du th\'eor\`eme 6
ci-dessous  fournit donc l'exemple cherch\'e.
\par
\bg
R\'esumons la construction suivante de [{\bf P2}], voir aussi
[{\bf P3}]:
\par
$ E_0, B$, et $ S $ sont des espaces de Banach, $ S $ est un 
sous-espace ferm\'e de $ B $, et $ i: S \hookrightarrow B $ 
est l'injection canonique. Soit $ u: S\rightarrow E_0 $ un
op\'erateur de norm $ \leq \eta \leq 1 $; alors  il existe
un espace de Banach $ E_1 $, un op\'erateur $ \wt u :
B \rightarrow E_1 $ et une injection isom\'etrique $
j : E_0 \hookrightarrow E_1 $, tels que 
$ \N{ \wt u } \leq 1 $ et $ \wt u \rond i = j \rond u $.
\bg
$$\matrix{B&\mapright{\wt u}&E_1\cr
       \mapupl{i}&{}&\mapupr{j}\cr
       \lower2.0ex\hbox{$S$}&\lower2.0ex\hbox{$
\mapright{u}$}&\lower2.0ex\hbox{$E_0$} } $$

\par
En effet, soit $ B_1 = B\oplus_1 E_0$, c'est $ B\times E_0 $ 
muni de la norme $ \N{(b,e)} = \N{b} + \N{e} $, soit le
sous-espace ferm\'e $ \Gamma(u) = \{(t,-u(t)): t\in S \}$ 
de $ B_1 $, et soit $ \pi : B_1 \rightarrow B_1/\Gamma(u)$  
 la surjection canonique. On pose $ E_1 = B_1/\Gamma(u)$,  
$ j(e) = \pi((0,e))$, et $ \wt u(b) = \pi((b,0))$; il est
alors  facile de v\'erifier que $ E_1 $, $j$ et $ \wt u$ 
ont les propri\'et\'es requises.
\par
\proclaim Th\'eor\`eme 3. 
Avec les m\^emes notations, on suppose que $(B,S) $ a $
{\cal P}(c) $, que  $ B$ est $H_1
(K,\sigma,I) $, $ E_0$ est $H_1
(K_{E_0},\sigma_0,I_0) $, et $ \eta \leq (1+c)/2 $. Alors,
l'espace $ E_1$, obtenu dans la construction pr\'ec\'edente,
 est un espace $H_1(K_{E_1},\sigma_1,I_1)$ avec
$$
I_0 \subset I_1,\qquad j \rond \sigma_0 = \sigma_1 \rond
s,\quad\hbox{et}\quad   K_{E_1} = \max (K_{E_0} , 2 c K).$$
\par
Preuve: 
Soit $ I_1 $ la r\'eunion disjointe de $ I $ et $ I_0 $. 
  Soit $ \sigma^\prime : \ell^1(I_1) \rightarrow B_1 $ l'application
d\'efinie par $ \sigma^\prime(x) = (\sigma(x_{\vert I}), 
\sigma_0(x_{\vert I_0}))  $, et $\sigma_1 = \pi \rond
\sigma^\prime : \ell^1(I_1)\rightarrow E_1 $. Il est
imm\'ediat de voir que $ \sigma_1 $ est une surjection
m\'etrique.
\bg
$$\matrix{B&\mapright{\wt u}&E_1&\mapleft{\sigma_1}&\ell^1(I_1)\cr
       \mapupl{i}&{}&\mapupr{j}&{}&\mapupr{s}\cr
       \lower2.0ex\hbox{$S$}&\lower2.0ex\hbox{$
\mapright{u}$}&\lower2.0ex\hbox{$E_0$}&\lower2.0ex\hbox{$
\mapleft{\sigma_0}$}&\lower2.0ex\hbox{$\ell^1(I_0)$} } $$

 \par
Soit $ f = \sum_0^m e^{ik(.)} a_k \in \wt H^1(E_1) $ de
norme  $ \N{f}_{\wt H^1(E_1)}  < 1 $; pour tout $ k $ on a 
 $ a_k = \pi ((y_k,\bar e_k)) $.
D'apr\`es $ {\cal P}(c) $, il existe $ 
g = \sum_0^{m_1} e^{ik(.)} x_k \in \wt H^1(B) $ telle que 
$$ \wt Q(g) = \sum_{k\geq 0} e^{ik(.)} Q(x_k) =
 \sum_{k\geq 0} e^{ik(.)} Q(y_k) $$
et
$$\N{g}_{\wt H^1(B)} \leq c \N{\sum_{k\geq 0}
e^{ik(.)} Q(x_k)}_{\wt H^1(B/S)} \eqno (2)$$
o\`u $ Q : B \to B/S $ est le quotient canonique.
\par
 On pose alors, pour $ k \geq 0 , e_k = \bar e_k - u(x_k -
y_k )$, d'o\`u  $ a_k = \pi ((x_k, e_k)) $, mais $ \N{f}_{\wt
H^1(E_1)} < 1 $,  donc il existe $ t : \tore \to S $
mesurable, telle que 
$$\int_\tore \N{\sum_{k=0}^{m_1} e^{ik\th} x_k + t(\th) }_{\wt
H^1(B)}\, dm(\th) +
\int_\tore \N{\sum_{k=0}^{m_1} e^{ik\th} e_k - u(t(\th)) }_{\wt
H^1(E_0)}\, dm(\th) < 1 \eqno (3) $$
notons $ \alpha, \beta $ respectivement la premi\`ere  et la 
seconde des int\'egrales pr\'ec\'edentes. En utilisant (2)
on a 
$$ \int_\tore \N{\sum_{k=0}^{m_1} e^{ik\th} x_k }_{\wt
H^1(B)}\, dm(\th) \leq c \alpha \eqno (4) $$
et
$$\int_\tore \N{ t(\th) }_{\wt H^1(B)}\, dm(\th)  \leq (1 +
c) \alpha \eqno (5) $$
et donc d'apr\`es (3) et (5) 
$$\int_\tore \N{\sum_{k=0}^{m_1} e^{ik\th} e_k  }_{\wt
H^1(E_0)}\, dm(\th) \leq \beta + \eta (1 + c )\alpha
\leq \beta +{ \alpha\over 2 }. \eqno (6)$$
\par
Il en r\'esulte que pour $
g = \sum_0^{m_1} e^{ik(.)} x_k \in \wt H^1(B) $ (resp. 
$g_0 = \sum_0^{m_1} e^{ik(.)} e_k \in \wt H^1(E_0) $),  il
existe $ h \in \wt H^1(\ell^1(I)) $ (resp. 
$ h_0 \in \wt H^1(\ell^1(I_0)) $)
telle que $ \wt
\sigma (h) = g $ et $ \N{h}_{\wt H^1(\ell^1(I))} \leq K
\N{g}_{\wt H^1(B)}  $ (resp.$ \wt
\sigma_0 (h_0) = g_0 $ et $ \N{h_0}_{\wt H^1(\ell^1(I_0))}
\leq K_{E_0} \N{g_0}_{\wt H^1(E_0)}  $).
\par
On d\'efinit alors,
$$ h_1(j,\th)=\left\{ \matrix{h(j,\th) &\hbox{ si $j \in I$}
\cr{}\cr
h_0(j,\th) &\hbox{ si $j \in I_0$}}\right .$$
\par
On a imm\'ediatement $ \wt \sigma_1(h_1) = f $ et
$$ \N{h_1}_{\wt H^1(\ell^1(I_1))} \leq K \N{g}_{\wt H^1(B)} 
+ K_{E_0} \N{g_0}_{\wt H^1(E_0)} \leq \max ( K_{E_0}, 2 c K)
$$
d'apr\`es (4), (6) et le fait que $ \alpha + \beta < 1 $.
\par
\bg
En raisonant par homog\'en\'eit\'e, on obtient alors le corollaire
suivant:
\proclaim Corollaire 4.
On suppose que $ (B,S) $ v\'erifie ${\cal P}(c)$, 
que  $ B$ est $H_1(K,\sigma,I)$, $ E_0$ est $H_1
(K_{E_0},\sigma_0,I_0) $, et $ u : S \to E_0 $. Alors il
existe $E_1 $, une injection isom\'etrique $ j: E_0
\hookrightarrow E_1 $ et un op\'erateur $\wt u : B \to E_0$ 
tels que
\item{$1.$} $E_1 \hbox{ est } H_1(K_{E_1},\sigma_1,I_1)$,
 avec $ K_1 = \max ( K_{E_0}, 2 c K)$), $I_0 \subset I_1$
et $ \sigma_1 \rond s = j \rond \sigma_0$.
\item{$2.$} $ \wt u_{\vert S } = j \rond u $, et $ \N{\wt u }
\leq 2 (1 + c) \N{u}$.
\par
\bg
En suivant Pisier dans [{\bf P2}], on en d\'eduit:
\proclaim Corollaire 5.
Soit $ \{ B_\lambda\}_{\lambda\in \Lambda} $ une famille
d'espaces de Banach, telle que, pour tout $\lambda\in
\Lambda$, l'espace
 $ B_\lambda$ est $H_1(K_{B_\lambda}, \sigma_\lambda,J_\lambda)$, et 
$ S_\lambda $ un sous-espace ferm\'e de $ B_\lambda $
tel que 
$ (B_\lambda, S_\lambda) $  ait $ {\cal P}(c_\lambda)$.
On suppose
$$
K = \sup \left\{ K_{B_\lambda} :  \lambda\in \Lambda\right
\} < + \infty      \qquad\hbox{ et }\qquad    c = 
\sup \left\{ c_\lambda :  \lambda\in \Lambda\right
\} < + \infty .$$
On consid\`ere, d'autre part, un espace de Banach $ E_0 \hbox{ v\'erifiant } H_1
(K_{E_0},\sigma_0,I_0) $ et une famille
d'op\'erateurs $ \{ u_\lambda : S_\lambda
\to E_0 \}_{ \lambda\in \Lambda } $.  Alors, Il existe $ E_1$, 
 une injection isom\'etrique $ j: E_0
\hookrightarrow E_1 $, tels que\hfill\break
\indent 1. $ E_1 \hbox{ est } H_1(K_{E_1},\sigma_1,I_1) $ 
(avec $ K_1 = \max ( K_{E_0}, 2 c K) $) , $ I_0 \subset I_1 $
et $ \sigma_1 \rond s = j \rond \sigma_0 $\hfill\break
\indent 2. Pour tout $ \lambda\in \Lambda $, il existe 
$ \wt u_\lambda : B_\lambda \to E_1$ qui v\'erifie
$\wt u_{\lambda\vert S_\lambda } = j \rond u_\lambda $,
et $ \N{\wt u_\lambda } \leq 2 (1 + c) \N{u_\lambda}$.
\par
\bg
Preuve: 
On se ram\`ene par homog\'en\'eit\'e \`a $ \forall \lambda\in 
\Lambda , \N{u_\lambda} \leq 1 $, et on pose 
$$
B = \ell^1( \{ B_\lambda\}_{\lambda\in \Lambda} ) = 
\left \{ (x_\lambda)_{\lambda\in \Lambda}  : x_\lambda \in
B_\lambda \hbox{  et } \sum_\Lambda \N{x_\lambda} < +
\infty \right \} $$
et
$$
 S = \ell^1( \{ S_\lambda\}_{\lambda\in \Lambda} ) = 
\left \{ (x_\lambda)_{\lambda\in \Lambda}  : x_\lambda \in
S_\lambda \hbox{  et } \sum_\Lambda \N{x_\lambda} < +
\infty \right \}. $$
On v\'erifie que $ (B,S) $ a la propri\'et\'e $ {\cal P}(c)
$, et que $ B \hbox{ est } H_1
(K,\Sigma,J) $ ( o\`u $ J $ est la r\'eunion disjointe de la
famille $ \{ J_\lambda \}_{\lambda\in \Lambda} $, et $
\Sigma $ est d\'efinie par $ \Sigma((x_\lambda)_{\lambda\in
\Lambda}) =(\sigma_\lambda(x_\lambda))_{\lambda\in
\Lambda}$).
On d\'efinit alors, $ u : S \to E_0 $ par 
$$ 
u((x_\lambda)_{\lambda\in\Lambda}) = \sum_\Lambda u_\lambda
(x_\lambda), $$
et on applique le corollaire 4.
\par
\bg
\proclaim Th\'eor\`eme 6. 
Il existe une constante  num\'erique  $ \kappa $, telle que,
pour tout espace $ E_0 $ qui est $ H_1(\mu, \sigma_0, I_0)
$ avec $ \mu \geq \kappa $, on peut trouver un espace 
$ E_1 $ qui est $ H_1(\mu, \sigma_1, I_1) $, et une
injection isom\'etrique $ j : E_0 \hookrightarrow E_1 $, 
telle que\hfill\break
\indent 1. $ I_0\subset I_1 $ et
 $ \sigma_1 \rond s = j \rond \sigma_0 $.\hfill\break
\indent 2.  Pour tout $ u \in E_0 \otimes E_0 $ on a
$$ \N{ j\otimes j(u) }_{E_1\hat\otimes E_1} \leq c(\mu) 
\N{ u}_{E_0\check \otimes E_0 }. $$
\par
\bg
Preuve : 
Pour $ p \in [1,+\infty[$, on note $ L_n^p $ l'espace $
\comp^n $ muni de la norme 
$$\N{x}_p = \left ( {1\over n}\sum_{k=1}^n
\abs{x_k}^p\right)^{1/p},$$
et on note $ i_n : L_n^2 \to L_n^1 $ l'application
identit\'e. On sait par une variante d'un th\'eor\`eme de 
Ka$\check{\hbox{s}}$in ([{\bf P3}]; Chapitre 7), que l'on
peut trouver une d\'ecomposition orthogonale; $ L_{3n}^2 = 
D_n^1 \oplus D_n^2 \oplus D_n^3$, avec $ \dim D_n^k =
n,(k=1, 2, 3)$ et une constante $ \delta > 0 $, telles que,
pour $ k=1, 2 $ 
$$
\forall  x \in  D_n^k \oplus D_n^3,     \hbox{ on a
}   \delta \N{ x }_2 \leq \N{i_{3n}(x)}_1 \leq \N{ x }_2 .$$
Pour $ k=1, 2$  on consid\`ere le quotient 
$$ Q_n^k :L_{3n}^1 \to B_n^k = L_{3n}^1/i_{3n}(D_n^k), $$ 
et on pose $ S_n^k = Q_n^k(D_n^3) $. 
\par
\noindent Remarquons que 
$$ B_n^k/S_n^k =  L_{3n}^1/i_{3n}(D_n^k \oplus D_n^3),$$
donc en utilisant le lemme ci-dessous (voir [{\bf K}]
pour une d\'emonstration utilisant un r\'esultat de [{\bf
BD}]),  on  voit facilement que les $ B_n^k $ sont 
$H^1$-projectifs, avec des constantes major\'ees
ind\'ependamment de $ n $, et que les couples $
(B_n^k,S_n^k) $ v\'erifient ${\cal P}(c) $ pour un $ c $
ind\'ependant  de $ n$.
\par
\proclaim Lemme. 
On suppose qu'il existe $\delta > 0 $, et des sous-espaces 
 $ Y_n \subset L_n^1 $ tels que $ \forall n ,\forall x \in
Y_n :    \delta \N{x}_2 \leq \N{x}_1$. Alors, il existe $ K $
telle que, si $ \sigma_n : L_n^1 \to L_n^1/Y_n^{} $ est
l'application  quotient, on a 
$$
\forall f \in \wt H^1(L_n^1/Y_n^{}),\quad  \exists g \in \wt
H^1(L_n^1)~ :~  \wt \sigma_n(g) = f  \qquad\hbox{ et }\qquad \N{g}_
{\wt H^1(L_n^1)} \leq K \N{f}_{\wt H^1(L_n^1/Y_n^{})}. $$
Soit $\{ u_{k,n,\ell}\}_{\ell \in \Lambda_n^k }$ la famille
des op\'erateurs de $ S_n^k $ dans $ E_0 $. On pose 
$$
\Lambda = \bigcup \left \{ (k,n) \times \Lambda_n^k :  k =
1, 2  \hbox{ et } n \geq  1 \right \}. $$
Si $ \lambda = (k,n,\ell) \in \Lambda $,  on d\'efinit 
 $$ (B_\lambda, S_\lambda,
u_\lambda) = (B_n^k,S_n^k,u_{k,n,\ell}) .$$
\par
Le corollaire 5 s'applique alors \`a cette famille, et on
trouve $  E_1 $ v\'erifiant le premier point du
th\'eor\`eme;  par contre le deuxi\`eme point est
d\'emontr\'e en utilisant le point 2. du corollaire 5  dans 
[{\bf P3}], p.141.
\par
\bg
\proclaim Th\'eor\`eme 7. 
Tout espace de Banach $H^1$-projectif $ E $ est
isom\'etriquement contenue dans un espace de Banach
$H^1$-projectif $ X $, tel que $ X \hot X = X \wcheck\otimes
X $.
\par
En effet, si $ E_0 = E $ est $ H^1(\mu,\sigma_0,I_0)$, alors 
en applicant le th\'eor\`eme 6 \`a $ E_0 $ puis \`a $ E_1 $, 
$E_2 $ \etc, on construit $ \{(E_n,j_n)\}_{n\geq 1} $ 
o\`u $ j_n : E_n \hookrightarrow E_{n+1} $ est une injection
isom\'etrique,$ E_n $ est $ H^1(\mu, \sigma_n, I_n)  $ avec
$ I_n \subset I_{n+1} $ et $ j_n \rond \sigma_n
=\sigma_{n+1} \rond s_n $;
\par
$$\matrix{E_0&\mapright{j_0}&E_1&\mapright{j_1}&\ldots&E_n&\mapright{j_n}&
E_{n+1}&\mapright{j_{n+1}}&\ldots\cr
\mapupr{\sigma_0}&{}&\mapupr{\sigma_1}&{}&{}&\mapupr{\sigma_n}&{}&
\mapupr{\sigma_{n+1}}&{}&{}\cr
\lower2.0ex\hbox{$\ell^1(I_0)$}&\lower2.0ex\hbox{$\mapright{s_0}$}&  
\lower2.0ex\hbox{$\ell^1(I_1)$}&\lower2.0ex\hbox{$\mapright{s_1}$}&
\lower2.0ex\hbox{$\ldots$}&
\lower2.0ex\hbox{$\ell^1(I_n)$}&\lower2.0ex\hbox{$\mapright{s_n}$}&  
\lower2.0ex\hbox{$\ell^1(I_{n+1})$}&\lower2.0ex\hbox{$\mapright{s_{n+1}}$}&
\lower2.0ex\hbox{$\ldots$}}$$
\noindent de plus
$$\forall u \in E_n\otimes E_n,     \hbox{ on a }   \N{
j_n\otimes _n (u) }_{E_n\hat\otimes E_n} \leq  c(\mu) 
\N{u}_{E_{n+1}\check\otimes E_{n+1}}.$$
\par
Prenons $ X $ la limite inductive de $ \{(E_n,j_n)\}_{n\geq
1} $ qui peut \^etre identifi\'ee \`a $ \adh{\cup E_n}$.
Notons $ I = \cup I_n $. Il est clair que $\cup
\ell^1(I_n) $ peut \^etre consid\'er\'e comme un
sous-espace dense dans $ \ell^1(I) $. On d\'efinit alors 
$\sigma : \ell^1(I) \to X $ par $ \sigma( \alpha) =
\sigma_n(\alpha) $ si $ \alpha \in \ell^1(I_n)$.
\par
En utilisant des arguments standards de densit\'e, on
d\'emontre que $ \wt \sigma $ est surjectif, et que 
$ X \hot X = X \wcheck\otimes X $. Ce qui d\'emontre le
th\'eor\`eme.
\par
\noindent{\it Remarques. }\par
\item{--} Il n'est pas difficile de modifier la construction 
pe\'ec\'edente pour trouver $ X $ v\'erifiant le
th\'eor\`eme 7, et tel que $ X $et son dual$ X^\ast $ soient
des espaces de cotype 2, v\'erifiant le th\'eor\`eme de
Grothendieck.
\item{--} Cet exemple montre aussi que le produit tensoriel
projectif de deux espaces Hardy-convexifiables(\cf [{\bf
X}]) n'est pas, en g\'en\'eral,  
Hardy-convexifiable.
\item{--} On peut aussi construire un espace de Banach $X$,
$H^1$-projectif,  ayant la propri\'et\'e d'approximation,
et tel que $ X \hot X $ ne soit pas $H^1$-projectif.
\par
\centerline{\sc R\'ef\'erences }
\par
\parindent40pt
\ref BD & J.  Bourgain and W.J. Davis
&  Martingales transforms and complex
uniform convexity, & Trans. Amer.
Math. Soc.& 294 (1986), 501--515.
\par
\ref HP & U. Haagerup and G. Pisier &
Factorization of analytic  functions 
with values in non-commutative
$L^1$-spaces, & Canad. J. Math. &41 (1989), 882-906.
\par
\ref K & O. Kouba & $H^1$-projective spaces& Quat. J. Math. Oxford (2)& 41(1990), 295-312.
\par
\ref P1 & G. Pisier & Factoriztion
of  operator valued analytic
functions, & Advances in Math.&93 No 1, (1992), 61-125.
\par
\ref P2 & G. Pisier & Counterexamples to a conjecture of
Grothendieck, & Acta Math.&151, (1983), 181--208.
\par
\ref P3 & G. Pisier && Factorization 
of  linear operators and geometry of
Banach spaces, & CBMS No 60, A.M.S.
Providence (1987).
\par
\ref X & Q. Xu & In\'egalit\'es pour 
les martingales de Hardy  et
renormage des espaces quasi-norm\'es, & C.R. Acad. Sci. Paris t. 307
S\'erie I, & (1988), 601--604.
\par
\end